\documentclass[12pt]{iopart}

\usepackage{amssymb}
\usepackage{graphicx}
\usepackage{iopams}  
\usepackage{xcolor}
\usepackage{tikz-cd}
\usepackage{ulem}

\newtheorem{theorem}{Theorem}
\newtheorem{proposition}{Proposition}

\newtheorem{lemma}{Lemma}
\newtheorem{corollary}{Corollary}

\def\reals{\blackboardbold{R}}

\usepackage{soul}

\def\reals{\mathbb{R}}

\def\phi{\varphi}

\def\bG_N{G_N}
\def\reals{\mathbb{R}}

\def\Om1{\Omega^{-1}}

\def\phi{\varphi}

\def\bG_N{G_N}
\def\param{\theta}

\def\catm{\mathcal{M}}
\def\phase{P}
\def\x{\mathbf{x}}
\def\y{\mathbf{y}}
\def\c{\mathbf{c}}

\def\->{\rightarrow}
\def\|->{\mapsto}

\def\C{\mathcal{C}}

\def\change{\color{red}}

\def\qed{$\blacksquare$} 
\def\params{P}
\def\phase{M}

\def\nhead{\vspace{4pt}\noindent}
\def\catmsing{\mathcal{S}}
\def\paramsing{\mathcal{B}}

\def\bC{\mathcal{C}}
\def\change{}

\begin{document}

\title[Existence of cusps]{Bistable boundary
conditions implying codimension 2 bifurcations}

\author{D.\ A.\ Rand}

\address{Mathematics Institute \& Zeeman Institute
for Systems Biology and Infectious
Disease Epidemiology Research,
University of Warwick,
Coventry CV4 7AL, UK}
\ead{d.a.rand@warwick.ac.uk}
\noindent\author{M.\ Saez}

\address{IQS, Universitat Ramon Llull, Via Augusta 390, 08017 Barcelona, Spain}
\ead{meritxell.saez@iqs.url.edu}
\vspace{10pt}
\begin{indented}
\item[]December 2024
\end{indented}

\begin{abstract}
We consider generic families $X_\param$
of smooth dynamical systems depending on
parameters $\param\in P$ where $P$ is a
2-dimensional simply connected domain
and assume that each $X_\param$ only has
a finite number of restpoints and
periodic orbits.  We prove that if over
the boundary of $P$ there is a S or Z
shaped bifurcation graph containing two
opposing fold bifurcation points while
over the rest of the boundary there are
no other bifurcation points, then, if there is no
fold-Hopf bifurcation in $P$ then there
is at least one cusp in the interior of
$P$.
\end{abstract}



%
%
%
%
%

\begin{figure}[h]
\centering
	\includegraphics[width=0.6\linewidth]{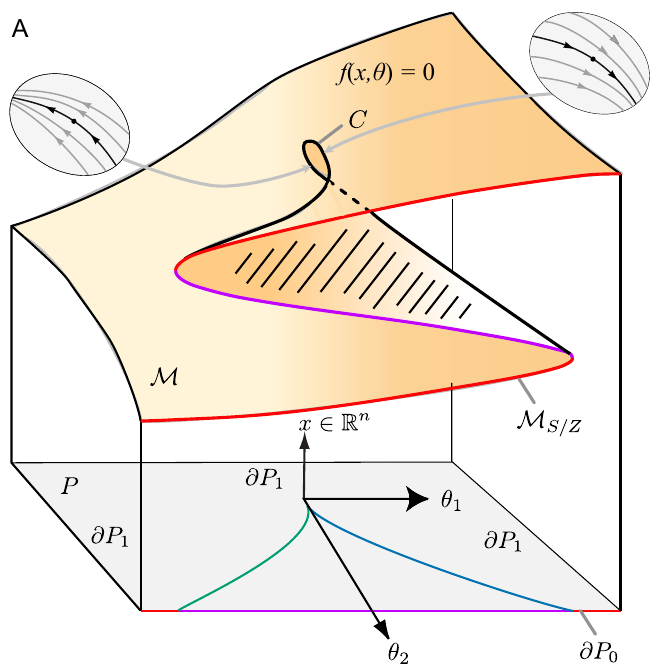}
	\caption{
	A. Over the boundary $\partial P$ of the parameter
	space there are no bifurcation points except the
	two folds in the S/Z curve over $\partial P$
	(red and purple curve). The folds are opposed
	in the sense defined in the text, a concept that
	formalises the notion of a S/Z curve. Theorem 1
	asserts that in this case there are an odd number of cusps
	in $P$. The figure shows the simplest case, where
	there is one. The coloured folded surface is the
	catastrophe manifold for this example. Note how the
	dynamics on the centre manifold of the fold point switches
	as the black fold curve goes through the cusp. 	
 } \label{fig:1}
\end{figure}

One of the most ubiquitous observations
in applied dynamical systems and many
areas of application is the S (or Z)
shaped bifurcation graph of the sort
shown in red and purple in Fig.\ \ref{fig:1}A.
The figure shows how the bifurcating restpoints
vary with the parameters.  
Such a 1-dimensional bifurcation graph can be found 
in almost any discussion of bistability
and is often discussed in a context where
there is more than one control parameter.
When the
parameter space is a 2-dimensional
simply connected domain $P$ it is often
the case that over its boundary there is
such a S or Z shaped bifurcation curve while
over the rest of the boundary there are no other bifurcation points.
It has
been assumed (e.g.\ in \cite{poston1978deducing}
and \cite{zeeman1977CT}) that when 
the system is gradient (as is 
assumed, for example, in catastrophe theory applications) under
reasonable and generic conditions
there must then be at least one cusp bifurcation in $P$. 
Indeed, this is a key assumption
in the huge number of applications of the cusp catastrophe
in the literature (see below) where one asserts the 
boundary conditions using 
data or otherwise and then deduces
the existence of a cusp catastrophe.
For example, it is a key assumption in
the formulation in \cite{ZeemanPrimary} of the Cooke-Zeeman
model for somitogenesis \cite{ZeemanCooke}
which  to quote a recent review \cite{McDaniel,francois2024} is still 
the most influential model in the field.
Similarly for Zeeman's still highly cited 
model of stock market dynamics \cite{zeeman1974unstable}.

In fact, we consider a more general
class of systems namely 2-parameter
families with only a finite number of
restpoints and periodic orbits. We show that, apart
from genericity, the conjectured result
is true without any conditions
(except finiteness) on the
dimensionality of the phase space or the
number of equilibria and periodic orbits
present.  The proof is necessary
because, although well-known local
bifurcation results classify the local structure of
the generic codimension 2 bifurcations
(\cite{GuckHolmes,arnold2013dynamical,kuznetsov2004}),
this still leaves the task of
showing that there must be such a point.

Besides its obvious general interest
from a theoretical point of view such a
result is very likely to be of use in
data science because there is much
interest in finding bifurcations in data
(e.g.\ for the cusp bifurcation in
data science, ecology/environment,
engineering, global warming, economics and developmental biology:
\cite{sguotti2020non,sguotti2020synergistic,
alamgir2019parametric,
shao2021method,
mostafa2020catastrophe,
huang2020risk,
zhang2020endpoint,
sun2020coupling,
saez2022statistically})
and codimension 2 points will be hard to
identify in noisy data.  On the other
hand, the attractor parts of S (or Z)
shaped bifurcation graphs should be much
easier to observe. Indeed, our interest
in this problem arose when using  
single cell data to analyse the
bifurcations that define the early
cellular decisions that stem cells make
in the early vertebrate embryo 
(e.g.\ \cite{saez2022statistically}).

Whether or not  this result is true for
gradient or gradient-like systems
was a key point of
contention during the controversy about
catastrophe theory in the 70s when it
was claimed in \cite{zahler1977} that it
is not true even under any reasonable
dynamical hypotheses (see also 
Smale's
review \cite{SmaleZeemanBook} of 
Zeeman's book \cite{zeeman1977CT}).  
It is therefore
remarkable that this claim has not been
clarified except in the special case
where the phase space is 1-dimensional
and the system is gradient
\cite{stewart1980catastrophe}.

This is a relatively rare example 
of  a global result about
local bifurcations in that it
is about systems rather than a local
result about germs. A condition on the boundary
implies certain bifurcations inside.
While catastrophe theory and local bifurcation theory
provide many powerful results which
are critical for applications,
many other applications need such an extension.
A lovely example of such a global result
proved in a completely different way
is Guckenheimer's condition in 
\cite{guckenheimer1977cusps} that implies
the existence of at least four cusps
for systems like Zeeman's Catastrophe Machine.

The ideas needed to prove our result
have much more general utility and
we will return to them in a later paper.
The approach is highly topological and a
key tool are the fold
approximating curves (defined in Appendix 2) that we construct
and the use of certain
bundles over curves 
in the catastrophe manifold
whose fibres are dynamical
objects such as centre manifolds.

{\change Figure \ref{fig:1} may help the reader
understand the reason we can probe the
parity of the number of cusps.  Each
point of the black fold curve on the
catastrophe manifold in Fig 1
corresponds to a fold bifurcation point
and this fold point has associated a
1-dimensional centre manifold.  As you
move along this fold curve from the
fold point at one end to the
other the tangent to this centre
manifold varies smoothly.  However, the
flow on the centre manifold switches direction at
the cusp point.  This is because the
folds at the end points are opposed in
terms of this flow.  If there were an
even number of cusps on the fold curve
then they would not be opposed.  By adding to this
curve the saddle part of the S/Z curve we get a circle
and it is circles like this and the corresponding
centre manifolds that we study.
For a general
MS family if we can find in $(x,\param )$
space a  closed curve $C$ of
restpoints which passes through all the
cusps and understand the topology of the
bundle whose fibre at each point
$(x,\param )$ in $C$ is the centre
manifold of the restpoint $x$, then we
can similarly get at the parity of the number of
cusps on the curve by counting the
switches.  We will find a curve $C$
where the bundle is a cyclinder so the
number of switches must be even.
However, one of these will correspond to
the opposite orientations of the fold
points in the S or Z curve and therefore
the number of switches corresponding to
cusps will be odd.}

\section{Main results}
We consider 2-parameter families of
smooth dynamical systems whose
non-wandering set consists of just a
finite number of equilibria and periodic
orbits.  Such a family consists of
smooth dynamical systems (flows)
depending smoothly on parameters which
vary in a region $\params$ of $\reals^2$ 
with a piecewise smooth boundary
$\partial \params$.  These families are
of the form $\dot{x}=X_\param
(x)=X(x,\param )$ where $x\in\phase$ and
$\param\in\params$ and $X$ is C$^r$ for 
$r\geq 3$.  It is well known
that for an open and dense subset of
such systems, away from bifurcation
sets, these systems are Morse-Smale and,
generically, the bifurcation set
consists of curves that are smooth
except at a finite number of points and
at these points a codimension-two
bifurcation occurs such as a cusp
bifurcation.  The bifurcation curve through a fold point
is smooth at all the fold points on it. The possible codimension 2 points
on such a curve are cusps, Bodganov-Takens (BT), and fold-Hopf (fH) points
(\cite{GuckHolmes,arnold2013dynamical,kuznetsov2004}).
The curve in $P$ continues smoothly through BT and 
FH points but has a cusp singularity in $P$ at cusp points.

We assume that the phase space $\phase$
is an $n$-dimensional disk (i.e.\
diffeomorphic to $\{ x\in\reals^n :
||x||<1\}$) and that the flow is always
inwardly transverse to its smooth, topologically
spherical boundary $\partial\phase$.  

An important part of our analysis is
a study of the structure of the
\emph{catastrophe manifold} $\catm$
of the parameterised family $X_\param$
which is defined by
\begin{equation*}
\catm =\{ (x,\param): x \mbox{ is a restpoint of } X_\param\} \subset \phase \times \params
\end{equation*}
and the associated map $\chi : \catm
\rightarrow \reals^2$ defined by the
projection $\x=(x,\param )\mapsto
\param$.  Generically, the origin is a
regular point of $X$ and then $\catm$ is
a 2-dimensional submanifold of
$\reals^{n}\times\reals^{2}$.  
We call a hyperbolic saddle with a
$k$-dimensional unstable manifold a $k$-saddle and
we will be particularly interested in the 1-saddles.
We call
the lift to $\catm$ via $\chi$ of folds,
fold curves, and other bifurcation
points in $P$ by the same name.  
Thus,
we may without distiction refer to a fold point $\x$ in
$\catm$ or a fold point $x$ in $P$.

Suppose that the parameter space is the
square $\params=\{ \param = (\param^1,\param^2):|\param^i | <1\}$
with boundary $\partial \params$
and consider the subset $\partial\params_0$ 
with boundary $\partial \params$
in  where
$\param^1=1$. We consider
the restpoints sitting over
$\partial\params_0$ i.e.\
$\catm_{S/Z} 
= \{ (x,\param ): \param \in \partial\params_0 
\mbox{ and $x$ is a restpoint of } X_\param \}$.
We say that this is a
$S/Z$-curve with fold points
$\x_1$ and $\x_2$
if it satisfies the following conditions. 
\begin{itemize}
\item[(i)] It is a smooth simply connected curve which contains just two fold
points $\x_1$ and $\x_2$ and 
$\catm_{S/Z}\setminus \{\x_1,\x_2\}$
has three connected components two of which consist of attractors
and the other consists of index 1 saddles.
We call the latter the \emph{saddle curve}
and denote it by $\gamma_\mathrm{sad}$.
\item[(ii)] The two fold points $\x_1$ and $\x_2$ are
\emph{opposed} in the following sense.
Put an orientation on
the 1-dimensional unstable manifold
of one of the saddles in $\gamma_\mathrm{sad}$
and extend this orientation continuously to 
all the points $\x$ in $\gamma_\mathrm{sad}$.
Then \emph{the folds are opposed} if
the fold orientation of one 
of the fold points $\x_1$ and $\x_2$
agrees
with that of the saddles close to it
while the other disagrees with those close to it.
Clearly, this does not depend on the choice
of the orientation of the center manifolds.
\end{itemize}
Finally, we suppose that the boundary of 
$\params$ contains no other bifurcation points.

Consider the set $\mathcal{S}_{\mathrm{sad}}$ 
of points in $\catm$ that correspond to generic
1-saddles (i.e.\ hyperbolic saddles with a 1-dimensional unstable
manifold). This is clearly an open subset of $\catm$.
Let $\catm_S$ denote the connected component of
$\mathcal{S}_{\mathrm{sad}}$ that contains $\gamma_\mathrm{sad}$.

{\change
\begin{theorem}\label{thm:one}
For a generic family $X_\param$, 
if there is a single $S/Z$ curve over
$\partial \params$ with fold points 
$\x_1=(x_1,\param_1) $ and $\x_2=(x_2,\param_2)$
and no other bifurcation points
over $\partial \params$ then either there is
at least one generic fold-Hopf bifurcation
in $\params$
or the total number of cusp points in
the boundary of $\catm_S$ is odd.
In particular, there must be a generic
codimension 2 bifurcation in $P$.
\end{theorem}
}

\vspace{2mm}\noindent{\bf Notes.} 
1. It is not necessary to assume a bound on
the total number of restpoints. \enskip
2. Although it is assumed that the two fold
points $\x_1$ and $\x_2$ are the only
bifurcations near the boundary, other
bifurcations can be allowed away from
the boundary. \enskip
3. There are examples of systems
satisfying the hypotheses of the theorem
with any positive odd number of 
cusps. These are provided by the $A_{2k}$
catastrophes \cite{arnold2013cat}. \enskip
4. The condition in (i) above ensures that at least one of
the cusps is a standard cusp i.e.\
one where two attractors interact with a
1-saddle between them.

A pretty much immediate corollary of
Theorem 1 is a similar result in the
context of elementary catastrophe theory.
Instead of a vectorfield,
we consider a $C^3$ smooth family of functions
$f_\param (x) =  f(x,\param)$ with state $x\in M$
depending upon parameters $\param\in P$
as above.

In this case the definition of a $S/Z$-curve
is entirely similar to that above for a parameterised 
family of gradient-like
dynamical systems except that instead of restpoints we treat 
critical points of $f_\param$ (i.e.\
the restpoints of the gradient vectorfield of $f_\param$)
and we define the two fold points to be opposed if
the following condition holds. 
Choose one of the saddles in the saddle curve
and put an orientation on
the 1-dimensional eigenspace $E(\param )$ of 
the nonnegative eigenvalue of
the Hessian of $f_\param$ at this saddle.
Extend this orientation continuously to 
all the corresponding eigenspaces $E(\param )$ 
of the saddles in the saddle curve
and extend it to the fold points. We say that
\emph{the folds are opposed} if
 near one fold point
$f_\param$ increases in the direction of the orientation
while near the other it decreases.

Finally, as above, we suppose that the boundary of 
$\params$ contains no other bifurcation points.

\begin{theorem}\label{thm:two}
Assuming genericity of the family $f$,
if there is a single $S/Z$ curve over
$\partial P$ with fold points 
$\x_1=(x_1,\param_1) $ and $\x_2=(x_2,\param_2)$
and no other bifurcation points
over $\partial \params$ then
the boundary of $\catm_S$ contains
an odd number of cusp points.
In particular, there must be a generic
cusp bifurcation in $P$.
\end{theorem}

\section{Proof of Theorems \ref{thm:one} and \ref{thm:two}.}

\nhead{\bf Preliminaries.}

In the generic case, the
set $\catmsing$ of fold, cusp, BT
and fH
points is a 1-dimensional smooth
submanifold of the 2-dimensional
submanifold $\catm$ and equals the set
of singularities of $\chi$ (Sect.\ 8 and Lemma 10.4
of \cite{kuznetsov2004}).  Its image $\paramsing
=\chi(\catmsing)$ is generically the set of local
bifurcation points in $\params$ that 
have a zero eigenvalue.  Moreover,
the tangent to
the centre manifold at a generic fold,
cusp or BT point varies smoothly as a
function of the point on these curves.

A curve of fold points in $\catm$ can be
smoothly extended through a a cusp, BT
point or a fH point on $\catm$ (e.g.\
Lemma 10.4, \cite{kuznetsov2004}).  It
follows that these curves in $\catm$,
that we call \emph{fold curves}, are
either open curves that leave $\catm$
through its boundary or are smooth
circles.  In the latter case we call
them \emph{fold circles}.

Cusps come in two forms: \emph{standard} and \emph{dual}.
At a standard cusp point two index $i$-saddles 
(e.g. attractors when $i=0$)
collide with a single $i+1$-saddle, while at a dual
one, two $i+1$- saddles (e.g.\ a saddle with a 1-dimensional
unstable manifold when $i =0$) collide with an
$i$-saddle.

An important difference between cusps
and BT points is the following.  At a
fold point $\x = (x,\param )$ in $\catm$,
the center manifold of $x$ and its
tangent space $\ell(\x)$ have a well
defined orientation defined by the
direction of the flow on the centre
manifold at points $x'\neq x$ near $x$.
It follows from the nonlinear conditions
that apply to a generic cusp or BT bifurcation 
(Sect.\ 8.7.2 and 8.7.4 of \cite{kuznetsov2004})
that, although the centre manifold tangent is
changing smoothly, as one passes
through a cusp on a fold curve
this orientation switches, while for a BT
bifurcation it does not.
It is this switch that allows us to determine the
parity of the number of cusps but not that of the
BT points.

We are particularly interested in the 1-saddles
and the network of attractors and  1-saddles 
has a nicely characterised structure (\cite{Smale60}). In 
particular, for each such 1-saddle $x$ the unstable
manifold $W^u(x)$ links the saddle to either one
or two attractors. The first case  is 
not relevant to the results we pursue and
so we will always assume that each such 1-saddle is
linked to two attractors. 
If we mention a saddle without specifying the dimensionality
then it is a 1-saddle.

\begin{figure}[h]
	\centering
	\includegraphics[width=0.8\linewidth]{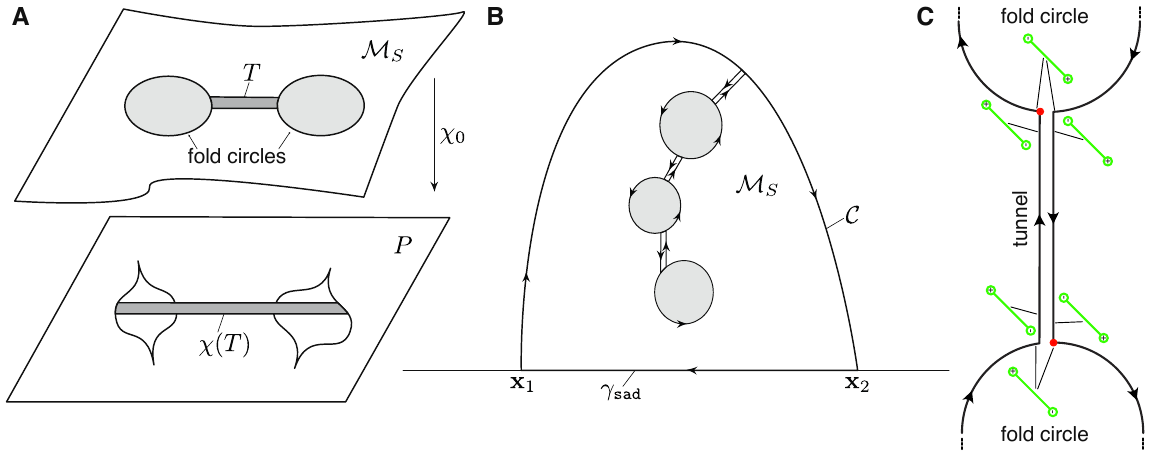}
	\caption{
	A. Schematic diagram of how a tunnel $T$ connects two fold
	circles in $\catm_S$.
	B. Schematic diagram of the curve
$\bar{\bC}$
	and how it is constructed using
	$\bC$ the fold circles and
	tunnels.  C.  A schematic
	illustration of why flips in
	orientation caused by the
	transitions from saddle to fold
	points in a tunnel occur in
	pairs.  The flips occur at the
	marked transitions from the
	tunnel.  Since there are no
	cusps near the tunnel ends the
	orientations at each of the
	tunnel ends agree.  Consequently,
	there is either zero or two
	flips in orientation as the two
	sides are traversed.
	}
	\label{fig:tunnels}
\end{figure}

Let $\catm_0$ be the connected component
of $\catm$ that contains $\catm_{S/Z}$
and let  $\chi_0$ denote $\chi_0 |_{\catm_0}$.
Since the two folds
$\x_i=(x_i,\param_i)$, $i=1,2$, in
$\catm_{S/Z}$ are generic there
is a fold curve crossing $\partial
\params_0$ transversally at the first
fold point $\param_1$.  This curve must
leave $\params$ and therefore must do
this at the second fold point $\param_2$
as there is no other fold point on the boundary.  
Let $\bC$
denote the lift via $\chi_0$ of this fold
curve to the catastrophe manifold
$\catm_0$.  This is a smooth curve.
The singular points of
$\chi_0$ in $\catm_0$
consist of $\bC$ and possibly
some disjoint fold circles. There
are no other open fold curves in $\catm_0$,
as otherwise there would be other fold
points on $\partial P$. 

Consider the
surface $\catm^*_\varepsilon$ 
given by  $\param_1 = x^3-x$, $0\leq \param_2 <\varepsilon$ 
in 3-dimensional $(x,\param_1,\param_2)$-space.
Let $\chi_0^*$ denote the restriction
to $\catm^*_\varepsilon$ of the projection
$(x,\param_1,\param_2 ) \rightarrow (\param_1,\param_2 )$.
In the argument below we will repeatedly use the
fact that there is a neighbourhood $N$ of
$\catm_{S/Z}$ in $\catm$ and diffeomorhisms 
$\phi$ of 
$N$ into $\catm^*_\varepsilon$ and $\eta$ of
$P$ 
into $\mu_1,\mu_2$ space such that
$\chi_0^*\circ \phi = \eta\circ\chi $.
We call this the local triviality
of $\chi_0$ near $\catm_{S/Z}$ and
it is explained further in Appendix 3.

Clearly, there is an annular
neighbourhood $N_{\partial P}$ of
$\chi_0^{-1} (\partial P)$ in $\catm_0$ such
that the only fold points in
$N_{\partial P}$ are two connected open
fold segments $S_1\ni \x_1$ and   $S_2\ni\x_2$ in $C\cap
N_{\partial P}$.  These arcs separate
$N_{\partial P}$ into a component
$N_{\partial P,S}$ consisting of saddles
and a component $N_{\partial P,A}$
consisting of attractors.

If $\y_1$ and $\y_2$ are two points in 
$\catm_0$ then $\y_1$ and $\y_2$  can be
connected by a smooth arc that enters
$N_{\partial P}$ and is transversal to
$S_1$ and $S_2$.  The parity of the
number of intersections is independent of 
such an arc. Consider the
equivalence relation: $\y_1\sim \y_2$
iff such curves connecting them have even parity
and note that the equivalence classes are
connected open sets and therefore
there are just two of them.
We deduce that $\bC$ separates $\catm_0$
into two connected components.  We label
these two components $\catm_S^\prime$
and $\catm_A^\prime$ according to
whether the arcs first enter
$N_{\partial P,S}$ or $N_{\partial
P,A}$.  Now consider the subset
$\catm_S$ (resp.\ $\catm_A$) of points
in $\catm_S^\prime$ (resp.\
$\catm_A^\prime$) that can be connected
to $N_{\partial P}$ as above by an arc
that does not contain any fold points.
Then all points in the subset 
have the same type and hence
are all saddles (in $\catm_S$) or all
attractors (in $\catm_A$).  It follows
that each fold circle in
$\catm_S^\prime$ (resp.\
$\catm_A^\prime$) separates
$\catm_S^\prime$ (resp.\
$\catm_A^\prime$ ) into two components
one of which contains $\catm_S$ (resp.\
$\catm_A$).  The other component is
called the interior of the fold circle.

For the proof of Theorem \ref{thm:one} we
focus on $\catm_S$.
Since there are no singularities of
$\chi_0$ in $\catm_S$, it is
orientable and has
no handles.
Therefore, by the
classification of surfaces (e.g.\ \cite{Hirsch2012}) they
are homeomorphic to the 2-sphere with a
number of (closed) disks removed.  
One of these holes corresponds to the
circle $\bC'$ formed by the union of $\bC$
and the saddle curve $\gamma_{\mathrm{sad}}$.
The others correspond to fold circles.

Consider a continuous parameterisation $\param(t)$, 
$0<t<t_\mathrm{end}$, of a closed curve in
$P$ that traverses the closed curve just once.
Then we say that that a lift $\Theta(t)$ of
$\param(t)$ to $\catm$ via $\chi$
is \emph{simple} if $\Theta(t_\mathrm{end})=\Theta(0)$.

In the following results we assume that there are
some fold circles in $\catm_S$. We deal with the
case where there are none in Proposition \ref{prop:cylinder}.
Moreover, we assume that there are no fH points on $\catm_S$ or its boundary as in that case the conclusion of Theorem 1 follows immediately.

{\change
The following lemma uses saddle approximating curves.
Given a fold curve in $P$ these are smooth curves that
are approximations of a bifurcation curve in $P$ that have a
lift to a smooth curve in $\catm$ via $\chi$ that approximates the
lift of the bifurcation curve. See Appendix \ref{sect:approx}
which contains a proof of their existence arbitrarily close 
to the bifurcation curve and its lift.

\begin{lemma}\label{lem:simple} 
Each of the fold circles in $\catm_S$
and its close saddle approximating curves
are the continuous lift via $\chi_0$ of a closed 
curve in $P$ and this lift is unique and simple.	
\end{lemma}

\nhead{\bf Note.} In general, there
are  examples of fold
circles where the fold circle has a
simple lift but no close approximation
curve does (e.g.\ 2-parameter families
close to the elliptic umbilic
catastrophe that have a bifurcation set
which is a closed curve with three cusps
on it).  The saddle approximating
curves surrounding its lift to $\catm$
double cover their images in $P$
so are not simple. 
It is the topological simplicity
of the catastrophe manifold near the
boundary that stops such behaviour in this lemma.

A $\varepsilon$-\emph{tunnel} in $\catm_S$ is a
pair of nonintersecting curves
$\x_i(t)$, $i=1,2$, that are C$^2$
$\varepsilon$-close and which project
under $\chi_0$ to a pair of
nonintersecting curves
(see Fig.\ \ref{fig:tunnels}).  The tunnel's
\emph{end segments} are
two short smooth segments
that join the end points so that the
tunnel and its end segments make up a
nonintersecting closed curve called 
the \emph{tunnel boundary}.  Note that if
a $\varepsilon$-tunnel in $\catm_S$ contains a point
$\x$ with $\x$ the only point
of $\catm_S$ in $\chi_0^{-1}(\chi_0(\x))$
then provided $\varepsilon >0$ is small enough
the lift of the tunnel boundary starting at $\chi_0(\x)$
via $\chi_0$ to $\catm_S$ is unique
and simple. Such a closed tunnel separates 
$\catm_S$ into two connected components
one of which meets the boundary of $\catm_S$.
The other component
is called  the \emph{interior} of the tunnel.

Given the topology of $\catm_S$
and the fact that $\chi_0$ is regular on it,
it follows that any two nonintersecting 
short smooth segments in 
$\catm_S$ can be connected
by a $\varepsilon$ tunnel in $\catm_S$ for 
any $\varepsilon > 0$ with the property that the
tunnel $T$ is the simple lift via $\chi$ 
of $\chi(T)$.

\nhead{\bf Proof of Lemma \ref{lem:simple}.}
Take a smooth saddle approximating curve
$\gamma_S$ in $P$ for the fold circle and let
$\Gamma_S$ be its lift to $\catm_S$ via
$\chi_0$.  This is a smooth curve in
$\catm_S$.  Connect $\Gamma_S$ to the
saddle curve $\gamma_{\mathrm{sad}}$ by the
$\varepsilon$-tunnel $T$ and let
$\gamma_1\subset\Gamma_S$ and
$\gamma_2\subset\gamma_{\mathrm{sad}}$
be the end segments of the tunnel.  Now
form a closed curve $C'$ in $\catm_S$
by deleting the end segments.

Take a circuit around the closed 
curve $\chi_0(C')$
starting at an endpoint of $\gamma_2$.
Then the end of the lift of the circuit
has to be the lift of the start point
because it is in $\gamma_{\mathrm{sad}}$
and there is at most one saddle above
each boundary point of the parameter
space.  But this means that the lift has
transversed all of $\Gamma_S$ except the
small segment $\Gamma_1$.  Since
$\Gamma_1$ can be made arbitrarily small
it follows that $\Gamma_S$ is the simple
lift of $\gamma_S$.
This shows that
each close approximating curve has a
simple lift.  Consequently, the same is
true of the fold circle.

The uniqueness of the lift follows from the
uniqueness of the lift of the tunnel.
\qed

\begin{corollary}\label{cor:annulus}
Each fold circle $C$ in $\catm_S$
is a boundary circle of an open annulus $A$	
in $\catm_S$ on which
$\chi_0$ is injective.
\end{corollary}

\nhead{\bf Proof.}
If $\x$ and $\y$ are points close to $C$
and $\chi_0(\x) =\chi_0(\y)$ take an
approximating curve $\gamma_S$ and
its lift $\Gamma_S$ with 
$\chi_0(\x) =\chi_0(\y)$ the 
start point of the lift. Since the lift is unique
it must pass through both $\x$ and $\y$. Since it is simple
this implies $\x = \y$.
\qed

It follows from Corollary \ref{cor:annulus} that
the end segments  of a tunnel can
also be in the boundary of
$\catm_S$ with the rest of the tunnel in 
$\catm_S$. This is because
then if $\x$ is a point in 
such a segment there are coordinates 
$(u,v)$ on a neighbourhood of $\x$ in which
the segment near $x$ is given by $u=0$. Thus
the curves given by $v=\pm 0.5 \varepsilon$ 
can be used for that end of the tunnel.

Let the fold circles in $\catm_S$ be
denote by $C_i$, $i=1,\ldots ,s$.

\begin{proposition}\label{prop:null_curve}
For $\varepsilon$ sufficiently small the
fold curves $C_i$ can be modified by
joining $\varepsilon$-tunnels between
$C_i$ and $C_{i+1}$ so that the
resulting curve $\tilde{\bC}$ is a
closed curve without self-intersections
that is one of the two connected
boundary curves of an annular region in
$\catm_S$ which has the curve $\bC'$ as
its other boundary component.
\end{proposition}

\nhead{\bf Proof.}
By the above discussion, $\bar{\catm_S}$
has the topology of a closed 2-dimensional disk
with $s$ disjoint open disks removed.

For each $i=1,\ldots ,s$, take two small disjoint
segments $\gamma_i$ and $\gamma_i^\prime$ in $C_i$ that meet
no cusps. 
Connect $\gamma_i$ to
$\gamma_{i+1}^\prime$ by a
$\varepsilon$-tunnel $T_i$ and delete
the segments in $\gamma_i$ and $\gamma_{i+1}^\prime$ between the
two endpoints of $T_i$.  This gives a
closed curve $T$. Given the topology of
$\bar{\catm_S}$ removing the interior of a tunnel
has the effect on the topology of reducing the
number of holes by one. Thus removing the interior of all tunnels
gives an annulus.
\qed

We construct a curve $\bar{\bC}$ by
using a $\varepsilon$-tunnel to join
$\bC$ to the part of the 
fold circle $C_s$ in 
the curve $\tilde{C}$
of Prop.\
\ref{prop:null_curve}.  The
ends of the tunnel should be 
chosen so
that they contain no cusp points and the
tunnel minus its end segments should lie
in the annulus given by Prop.\
\ref{prop:null_curve} (see Fig.\ \ref{fig:tunnels}).  Clearly
the part of $\catm_S$ not in any of
the interiors of the tunnels 
used in this construction is a 
topological disk. 
It follows that
this curve $\bar{\bC}$ which 
is its boundary is homotopically
null since a topological disk is simply connected.

All points in $\bar{\bC}$ apart from
a finite number of either
cusp or BT points are fold points or saddle points.
We form a bundle $\mathcal{B}_{\bar{\bC}}$
	over $\bar{\bC}$
whose fibre $\ell (\x)$ at $\x=(x,\param) $
is the tangent space to the unstable manifold
if $x$ is a saddle, and the tangent to the centre manifold
at $x$ otherwise. The fibres have continuous and piecewise smooth
dependence, not smooth only at the points where $\bar{\bC}$
is not smooth.

If $\gamma$ is any smooth curve on $\catm_S$
then we can define a similar bundle 
$\mathcal{B}_{\gamma}$
but where the
fibres over $\x=(x,\param) \in \gamma$
are the tangent space to the unstable manifold of $x$.
If such a curve is C$^2$ close to $\bar{\bC}$
at the smooth points of $\bar{\bC}$ then
the topological type of $\mathcal{B}_{\bar{\bC}}$
and $\mathcal{B}_{\gamma}$ are the same.

\begin{proposition}\label{prop:cylinder}
If there are fold circles in $\catm_S$
let $\bar{\bC}$ be as in Prop.\ \ref{prop:null_curve}.
Otherwise let $\bar{\bC}$ be $\bC'$.
Then, $\mathcal{B}_{\bar{\bC}}$ is trivial 
(i.e.\ has the topology of a cylinder).
\end{proposition}

\nhead{\bf Proof.}
Take a C$^2$ curve $C_0$ in $\catm_S$ that
C$^2$-approximates the smooth parts of 
$\bar{\bC}$ or $\bC$ if there are no fold circles
in $\catm_S$. Then  $C_0$ can be smoothly homotoped
within $\catm_S$ to an arbitrarily small
smooth circle $C_1$ and during this homotopy each 
closed curve has a saddle bundle with the
same topology as that of $C_0$.
But if $C_1$ is sufficiently small
this bundle must be trivial as if this is
sufficiently small, then, for all $\x, \x'\in C_1$,
$\ell(\x)$ and $\ell(\x')$ are never normal.
Thus, $\mathcal{B}_{C_1}$
and hence also $\mathcal{B}_{\bar{\bC}}$ must be trivial.
\qed

\nhead{\bf Proof of Theorem \ref{thm:one}.}

We follow a path around the curve
$\bar{\bC}$ of Prop.\ \ref{prop:cylinder} that starts
at $\x_1$
and continues firstly around $\bC$ and then the
rest of $\bar{\bC}$ until returning to $\x_1$
(see Fig.\ \ref{fig:tunnels}).
As we traverse the curve, we
will put a direction on each fibre $\ell (\x)$
encountered.  Whenever the base point
$\x=(x,\param)$ on $\bar{\bC}$
is a fold point we put the fold
direction on the fibre.  This is 
given by the direction of the flow near
$x$ on the centre manifold of $x$.
Where the curve
transitions from a fold point to a
saddle point we continue the direction
so that it is locally constant.  Where
the curve transitions from a saddle
point to a fold point we change it to
the new fold orientation.  This
direction field will have a
discontinuity at each cusp but not at a BT point.
It may also have a discontinuity as the
curve exits a tunnel but each tunnel is
exited twice and if there is a
discontinuity at one exit then there is one
at the other.  Consequently, the parity
of the total number $N$ of cusps
points equals that of all such above
discontinuities.  There is one other
discontinuity and this occurs as the
base point arrives back at $\x_1$.  This
is because the folds at $\x_1$ and
$\x_2$ are opposed.  There is no
discontinuity at $\x_2$ as this is a
fold to saddle transition and so the
direction put on fibres over base points
in $\gamma_{\mathrm{sad}}$ is that of
$\x_2$ which is the opposite of that at
$\x_1$.  It follows that the parity of
$N$ and the total number of
discontinuities are different.  But
since the discontinuities are discrete
and since $\mathcal{B}_{\bar{\bC}}$ is
trivial (i.e.\ a cylinder) the total
number is even.  Consequently, $N$ is odd,
proving the theorem
\qed

\nhead{\bf Proof of Theorem 2.}
Take any gradient system with parameterised potential $f$.
Such a system has no BT or fH bifurcations.
Then the parameterised gradient system satisfies the hypotheses 
of Theorem 1. For this we 
need to check the steps where we used genericity.
This follows since the catastrophe manifold for the
critical points of $f_\param$ and
the fold curves on it are the same as those for the
restpoints of the gradient system.
\qed
} 

\section*{Appendix 1: Center manifolds and the center manifold bundle.}\label{sect:center_manifolds}

\subsection*{Center manifolds and smoothness}
For relevant
information about center manifolds see
\cite{hirsch2006invariant} Sect.\ 5A.  In
particular note that by Theorem 5A.3 of
\cite{hirsch2006invariant}, if $W^c$ is a
center manifold through a restpoint $x$
and $W$ is a backward invariant set
containing $x$ then, near $x$, $W$ is
contained in $W^c$.  Thus, for example, if the
unstable manifold of a saddle is
asymptotic to a fold point, then close
to the fold point it is in the center
manifold.
Center manifolds are not necessarily
unique but their tangent space is.  We
will use this fact below.  

We now consider what we call
\emph{pseudo-hyperbolic} restpoints
$x$.  At such restpoints $x$ there is $a>b>0$
such that the Jacobian of the
vector field at $x$ has eigenvalues $\lambda$
that either have their real part $\leq
-a$ or $\geq -b$.  
Pseudo-hyperbolic index 1 saddles and
attractors have 1-dimensional center
manifolds $W^c(x)$ that vary smoothly
with parameters (Sect.\ 5
\cite{hirsch2006invariant}, especially
Theorems 5.1, 5.5 and 5A.1).  If $\varphi^{t}$
is the flow, this
manifold is characterised by the fact
that $z\in W^c(x)\iff ||
\varphi^{-t}(z)-x||/e^{ct} \rightarrow
0$ as $t\rightarrow \infty$ for any $c$
with $a>c>b$.  There is a complementary
submanifold $W^{ss}(x)$ transversal to 
$W^c(x)$ at $x$ characterised by
$z\in W^{ss}(x)\iff ||
\varphi^t(z)-x||/e^{-ct} \rightarrow 0$ as
$t\rightarrow \infty$ for such a $c$.
This we call the \emph{strong stable
	manifold}.  Note that our use of the
term \emph{center manifold} is a little
more general than usual as in that case
one commonly takes only $b=0$.

Index 1 saddles are always
pseudo-hyperbolic and attractors are if
they are close to having a fold
bifurcation.  For an index 1 saddle,
part of the unstable manifold containing
the saddle can be taken for a center
manifold.\label{para:pseudo-hyperbolic}

According to Theorem 5.1 of
\cite{hirsch2006invariant}, $W^c(x)$ has
$C^r$ dependence upon parameters
provided $e^{jb-a}<1$ for $1\leq j\leq
r$.  Thus the center manifold for
saddles always is smooth and that for
attractors is smooth provided they are close
enough to having a fold bifurcation.
The later point is true because
the closer an attractor is to being
a fold, the closer one can take $b$
to zero.

\subsection*{Approximations and the CM bundle}

Suppose we have a nonsingular curve
$\gamma(t)$, $0<t<T$ in either $P$ or
$\catm$ together with a tubular
neighbourhood $N$ of $\gamma$ and
consider another $C^r$ curve
$\tilde{\gamma}$ that passes through
$N$.  By definition of a tubular
neighbourhood there is a retraction
$\pi: N\->\gamma$ making $(\pi, N,\gamma
)$ a vector bundle whose zero section is
the inclusion $\gamma\-> N$.  We say
that $\tilde{\gamma}$ is
$\varepsilon-C^r$-close to $\gamma$ in
$N$ if the absolute value of the
derivatives of $\tilde{\gamma}(t)$ wrt
$t$ of order $0,\ldots ,r$ are within distance
$\varepsilon$ of those of
$\pi(\tilde{\gamma})(t)$.

When $\gamma\subset\catm$, we say that
$\gamma$ has center manifolds if at each
point $\x=(x,\param)$ of $\gamma$ there
is a $C^1$ center manifold at $x$ in $M$
and these vary $C^2$-smoothly with
$\x\in\gamma$.  We shall be especially
interested in the line bundle
$\mathcal{B}_\gamma$ over such a curve
$\gamma$ whose fibre at $\x=(x,\param)$
is the tangent space $\ell (\x)$ to the
center manifold at $x$.  By the above
discussion, if $\tilde{\gamma}$ is a
curve in $\catm$ that is sufficiently
$\varepsilon-C^2$-close to $\gamma$ then
the center manifolds for the restpoints
at $\tilde{\gamma}(t)$ and
$\pi(\tilde{\gamma}(t))$ vary in a $C^2$
fashion and their difference $d(\ell
(\pi(\gamma(t)))), \ell (\gamma(t))$ is
$O(\varepsilon)$ with the constant of
proportionality independent of $t$ if
the curve is compact.  Here $d(\ell
(\x), \ell (\x')) = \min ||e-e'||$ where
the minimum is over all unit vectors
$e\in \ell (\x)$, $e'\in \ell (\x')$.
Therefore, we have the following lemma.

\begin{lemma}
	If
$\gamma$ and $\tilde{\gamma}$ are closed
curves as above and $\varepsilon > 0$ is
sufficiently small,
 $\mathcal{B}_\gamma$ and
 $\mathcal{B}_{\gamma'}$
are both trivial bundles or they are
both topologically M\"obius bands.
\end{lemma}

\section*{Appendix 2: Approximating curves.}\label{sect:approx}
In the proof of Theorem 1 we only use
the saddle approximating curve.  The
existence of arbitrarily close attractor
approximating curves is also of use in other problems
but not needed here.  Indeed, near a BT bifurcation, there
are points where the attractor has a
complex conjugate pair of eigenvalues
and at these points there is no
1-dimensional centre manifold.  On the
other hand, a saddle approximating curve
behaves nicely at BT points and the 1-dimensional
unstable manifold of the saddle has
smooth variation.  However, this is not the case
for the fH bifurcation as along a curve close to the fold curve
the dimension of the unstable manifold of the saddle changes.

For the rest of this section by a fold curve we mean one
that involves a 1-saddle as this is all we are concerned about in the
proof of Theorem 1.

\begin{figure}[h]
	\centering
	\includegraphics[width=0.8\linewidth]{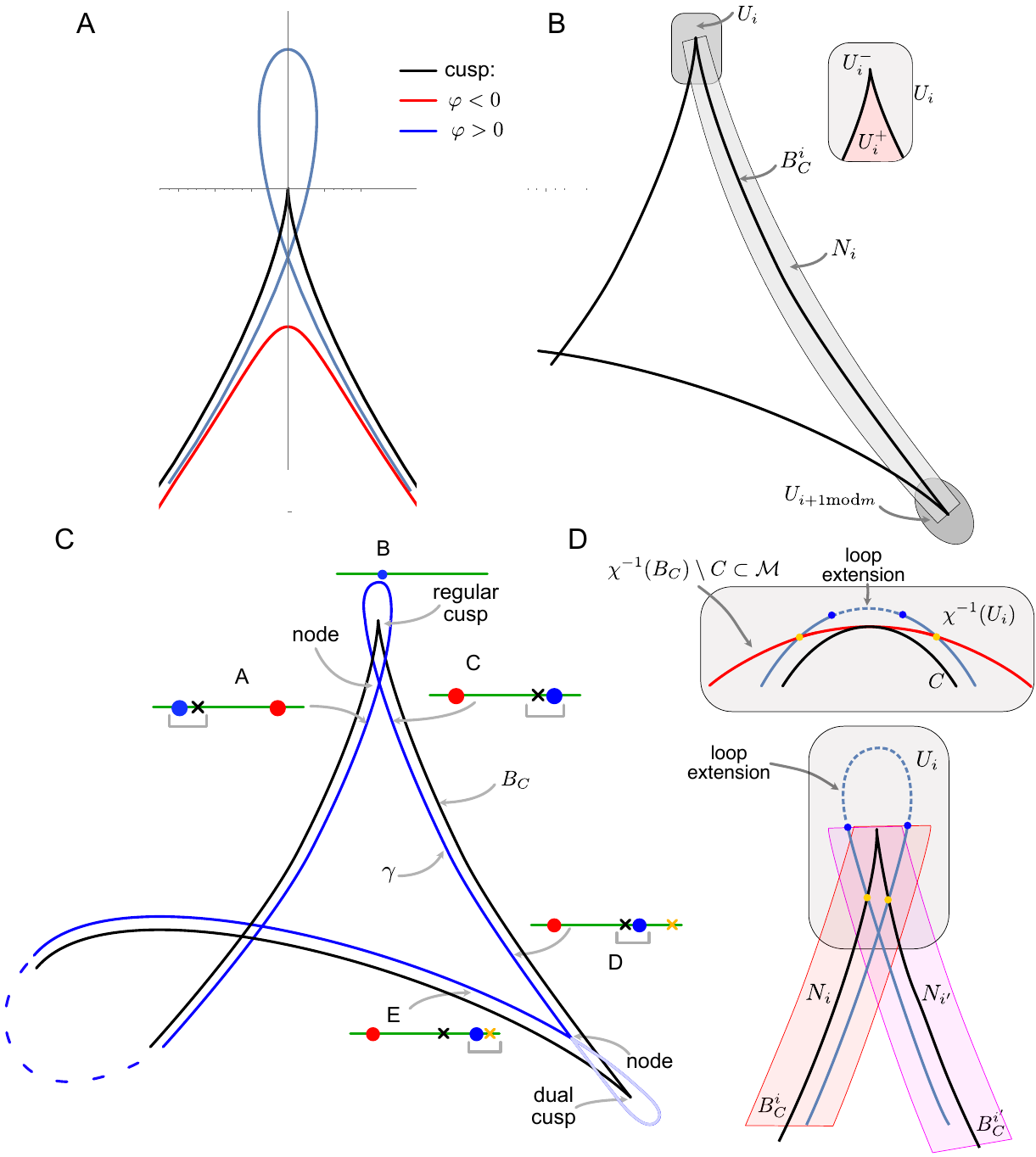}
	\caption{
	Schematic of neighbourhoods and approximating curves.
	A. Looping curve (red) and nudging curve (blue).
	B. Disposition of the neighbourhoods $U_i$,
	$U_{i \mathtt{ mod } m}$ and $N_i$.
	C. Showing how an approximating curve
	for an attractor 
	 goes around, first, around a standard  cusp
	 and then a dual cusp.
	 D. Explaining the loop extension.
	}
	\label{fig:3}
\end{figure}

The pair of equations
	$u = x^2$, $v = y,$
are the normal form for $\chi$ near a fold; (see
Theorem 15A, \cite{whitney1955singularities}).  
Therefore, if 
$\chi : \catm \rightarrow \params$ is
the mapping under consideration and
$\x\in\catm$ is a fold point,  there
is a neighbourhood $U$ of $\x$ in
$\catm$ and $V$ of $\chi(\x)$ in
$\params$ and a smooth curve $\gamma$ in
$V$ such that the lift via $\chi$ of
$\gamma$ to $U$ is a smooth curve
$\Gamma$ that is the set of fold points
in $U$.  The curves $\gamma$ and
$\Gamma$ separate $U$ and $V$
respectively into two connected
components, and one of the two components of $V$ does not
intersect $\chi (U)$. 

The equations $u = x^2$, $v = y,$ are also the
normal form for $\chi$ near a BT point on $\catm$.
This is because the normal form
for this codimension 2 bifurcation
is given by 
$ \dot{x}_1 = x_2,
\dot{x}_2 = \beta_1 + \beta_2 x_1 + x_1^2 +s x_1 x_2 +\mathrm{O}(||x||^3)$.

\def\u{\param_1}
\def\v{\param_2}
A normal form for the catastrophe manifold for
the standard cusp bifurcation is given by
the equation 
$x^3-\u x-\v=0$.  Therefore, the
map $(\u,x)\mapsto (x,\u ,\v=-\u x+x^3)$ from
$\reals^2$ to $\catm$ parameterises
$\catm$ in terms of $x$ and $\u$.  Thus,
in this parameterisation $\chi$ is given
by $(\u,x)\mapsto (\u,\v=-\u x+x^3)$ and this
is singular when $\u=3x^2$ which defines
a smooth curve $C$ in $\catm$.  The
bifurcation set $B_C$ is its image under
$\chi$, which is the set of points given
by $\u=3x^2,\v=- 2x^3$ i.e.\ $4\u^3=27\v^2$.
The dual cusp ($+$) case is entirely analogous.

Any curve in $\catm$ that is $C^r$-close
to $C$, $r>1$, is of the form
$\param_1=3x^2+\phi(x)$ and the image
under $\chi$ therefore has the
parametric form 
\begin{equation}\label{eqn:pf}
\param_1 =3x^2+\phi(x),\quad
\param_2 =-2x^3-x\phi(x).	
\end{equation}
Consequently, if
$\phi(x)$ is of constant sign the form
of the image curves are as shown in
Fig.\ \ref{fig:3}A.  In particular, if
$\phi >0$ then this curve is smooth and
loops around the cusp (blue curve in
Fig.\ \ref{fig:3}A) and if $\phi <0$ the
curve has no self-intersections and
stays inside the cusp (red curve).
We call these respectively
\emph{cusp looping curves} and
\emph{cusp nudging curves} for the cusp.
Conversely, any curve with the
parametric form (\ref{eqn:pf})
lifts via $\chi$ to a curve that
is $C^2$-close to $C$, if $\phi$
is $C^2$-small and furthermore,
in $\catm$, it lies to one
side of $C$. 

\begin{lemma}\label{lem:loop&nudge}
If $C$ is a generic fold curve,
$N$ a tubular neighbourhood of $C$,
$\x$ is a cusp point on $C$
and $\varepsilon>0$ then, in some neighbourhood
of $\chi(\x)$ there 
is a cusp looping curve $\gamma_\ell$ and
a cusp nudging curve $\gamma_n$
for $\x$ that lift via $\chi$ to
curves $\Gamma_\ell$ and $\Gamma_n$
that are $\varepsilon-\C^2$-close
to $C$ in $N$. If the cusp is
standard, $\Gamma_\ell$ will be a 
curve of attractors and $\Gamma_n$
a curve of saddles, and vice-versa
for a dual cusp.
\end{lemma}

\noindent{\bf Proof.}
We change coordinates to put the cusp in
normal form as above and then the result follows from the
discussion above.
\qed

The local structure of codimension 2 bifurcation
points is described in 
\cite{GuckHolmes}, \cite{arnold2013dynamical} and \cite{kuznetsov2004}

\begin{theorem}
For a generic MS family,
given $\varepsilon >0$ and a 1-saddle fold curve
$C$ in $\catm$ not containing any fH points there is a  $C^2$-
$\varepsilon$-approximating curve
$\gamma_S$  in $P$ with
the following property: $\gamma_S$
lifts via $\chi$ to
a $C^2$ curve $\Gamma_S$ (resp.\
$\Gamma_A$) of 1-saddles in $\catm$ that is
$\varepsilon-C^2$-close to $C$.  If $C$
is a fold circle then the $\gamma_S$ 
can be taken to be a closed
curve.  In this case taking the lift
via $\chi$ of the curve traversed $r$
times produces a closed curve in $\catm$
that is $\varepsilon^\prime-C^2$-close to $C$
where $\varepsilon^\prime$ is $O(\varepsilon)$.
If the fold curve contains BT points but not fH
points then such saddle approximating curves exist.
\end{theorem}

\noindent{\bf Proof.}
The key part of the proof is contained
in Lemma \ref{lem:loop&nudge}.
$B_C=\chi (C)$ will contain a possibly
empty set of finitely many cusps $\c_1
,\ldots , \c_m$, $\c_i=(c_i ,\param_i)$.
This labelling can be chosen so that
there are no cusps on $C$ between $\c_i$
and $\c_{i'}$ where $i'$ denotes $i+1$
if $C$ is open and
${i+1\mbox{\,mod}\,m}$ if $C$ is a fold
circle.  Let $C^i$ denote this segment
and $B_C^i$ denote $\chi (C^i)$.  Let
$N$ be a thin tubular neighbourhood of
$C$ in $\catm$ and $N_i$ be a thin
tubular neighbourhood of $B_C^i$ that
satisfies $\chi^{-1}(N_i)\subset N$.

We consider the arc $C^i$
of fold and BH points between $\c_i$ and
$\c_i'$.  By Lemma \ref{lem:loop&nudge} we can find a
cusp looping curve around each dual
cusp and a nudging curve at each standard
cusp which are in $N_i$ and which have
lifts that are sufficiently close to
$C$.  It is then straightforward to
join these by a curve inside each $N_i$
that lifts to a curve that is 
$C^2$-close to $C$ because,
as noted above, at
both fold and BT points, $\chi$ has 
the structure of a fold. In this way we
construct the curve. 
To see that traversing the curve multiple times
eventually gives a closed curve keep repeating the
above process each time taking for the
start point of the lift to $\catm$ the endpoint of the
previous lift.
\qed

\section*{Appendix 3. Local triviality
of $\chi$ near $\catm_{S/Z}$} 

Consider a system with state variable $x$ and
parameters $\mu_1$ and $\mu_2$. We consider the
surface $\catm^*_\varepsilon$ 
given by  $\mu_1 = x^3-x$, $0\leq \mu_2 <\varepsilon$ 
in $(x,\mu_1,\mu_2)$-space.
Let $\chi^*$ denote the restriction
to $\catm^*_\varepsilon$ of the projection
$(x,\mu_1,\mu_2 ) \rightarrow (\mu_1,\mu_2 )$.

\begin{lemma}
There is a diffeomorphism $\phi$
from a neighbourhood of
 	$\catm_{S/Z}$ in $\catm$
 	to $\catm^*_\varepsilon$
 	and a diffeomorphism $\eta$
 	between the two parameter spaces such that
 	$\chi^*\circ \phi = \eta\circ\chi$.
\end{lemma}

\def\U{\mathcal{U}}
\noindent{\bf Proof.}
Consider a thin tubular neighbourhood
$N$ of $\catm_{S/Z}$ in $\catm$.  Then
provided $N$ is thin enough, $N\setminus
\C$ has three connected components which
are discs. Moreover, there are two
neighbourhoods $\U_1$ and $\U_2$
respectively containing the two
connected components of $N\cap \C$
such that on $\U_i$, $i=1,2$, $\chi$ has
the normal form $(u_i,v_i)=(\pm
x_i^2,y_i)$ in some coordinate system
$(x_i,y_i)$.  The choice of sign will be different at the two fold curves in $N\setminus
\C$ since the fold points $\x_1$ and $\x_2$
are assumed to be opposed.

Let $D$ be the range $\chi
(N)$ together with a smooth structure
compatible with the
two sets of coordinates $(u_i,v_i)$,
$i=1,2$.  Then there is a diffeomorphism
of $\U_1\cup \U_2$  into two
neighbourhoods of the fold curves in
$\catm^*_\varepsilon$ such that the
diagram below commutes on $\U_1\cup
\U_2$.  

\begin{tikzcd}
\U_1\cup \U_2  \subset  N \arrow[r, "\phi"] \arrow[d, "\chi"]
& \catm^*_\varepsilon \arrow[d, "\chi^*"] \\
 P \arrow[r, "\eta"]
&  (\param_1,\param_2)-space
\end{tikzcd}

\noindent
Now we can extend the
diffeomorphims to $\catm_0$ and $D$
using the fact that, outside of
$\U_1$ and $\U_2$, the restriction of
$\chi$ to any one of the connected
components of $N\setminus \C$ is
injective.
\qed

\section*{References}
\bibliographystyle{ieeetr}

\end{document}